\documentclass[12pt,leqno]{amsart}
\pdfoutput=1
\usepackage{amsmath, amssymb, amsthm}
\usepackage[svgnames]{xcolor}
\usepackage{tikz}
\usetikzlibrary{arrows,cd}
\tikzset{>=latex}
\usepackage[shortalphabetic]{amsrefs}
\usepackage[T1]{fontenc}
\usepackage{mathptmx}
\usepackage{microtype}
\usepackage[scaled=0.75]{beramono}
\usepackage[colorlinks = true,
  linkcolor  = DarkBlue,
  urlcolor   = DarkRed,
  citecolor  = DarkGreen]{hyperref}
\usepackage{caption}
\usepackage{subcaption}
\usepackage[centering, includeheadfoot, hmargin=0.95in, tmargin=0.8in,
  bmargin=0.8in, headheight=6pt]{geometry}
\usepackage{listings}
\usepackage{stix}

\newcommand{\colequal}{\ensuremath{:\!=}}
\DeclareUnicodeCharacter{0393}{\Gamma}
\DeclareUnicodeCharacter{0394}{\Delta}
\DeclareUnicodeCharacter{22C8}{\bowtie}
\DeclareUnicodeCharacter{29D3}{\fbowtie}
\DeclareUnicodeCharacter{29D1}{\lfbowtie}
\DeclareUnicodeCharacter{2119}{\mathbb{P}}

\begin{document}

\vspace*{-4.5em}

\title[Simplicial Complexes]{Simplicial complexes in Macaulay2}

\author[B.~Hersey]{Ben Hersey}
\author[G.G.~Smith]{Gregory G.{} Smith} 
\author[A.~Zotine]{Alexandre Zotine}

\address{Department of Mathematics \& Statistics, Concordia
  University, Montréal, QB, H3G 1M8, Canada
  {\normalfont\texttt{benjamin.hersey@concordia.ca}}.
}

\address{Department of Mathematics \& Statistics, Queen's
  University, Kingston, ON, K7L 3N6, Canada
  {\normalfont\texttt{ggsmith@mast.queensu.ca}},
  {\normalfont\texttt{18az45@queensu.ca}}.
}

\thanks{
  \indent
  \href{http://www2.macaulay2.com/Macaulay2/doc/Macaulay2-1.20/share/doc/Macaulay2/SimplicialComplexes/html/index.html}{\texttt{SimplicialComplexes} version 2.0}}

\begin{abstract}
  We highlight some features of the \emph{SimplicialComplexes} package in
  \emph{Macaulay2}.
\end{abstract}

\maketitle

\vspace{-0.5em}

\noindent
This updated version of the \emph{SimplicialComplexes} package in
\emph{Macaulay2}~\cite{M2}, originally developed by Sorin Popescu, Gregory
G. Smith, and Mike Stillman, adds constructors for many classic examples,
implements a new data type for simplicial maps, and incorporates many
improvements to the methods and documentation.  Emphasizing combinatorial and
algebraic applications, the primary data type encodes an abstract simplicial
complex---a family of subsets of a finite set that is closed under taking
subsets.  These simplicial complexes are the combinatorial counterpart to
their geometric realizations formed from points, line segments, filled-in
triangles, solid tetrahedra, and their higher-dimensional analogues in some
Euclidean space.  The subsets in a simplicial complex are called faces. The
faces having cardinality~$1$ are its vertices and the maximal faces (ordered
by inclusion) are its facets.  The dimension of a simplicial complex is one
less than the maximum cardinality of its faces. Following the combinatorial
conventions, every nonempty simplicial complex has the empty set as a face.

In this package, a simplicial complex is represented by its Stanley--Reisner
ideal. The vertices are identified with a subset of the variables in a
polynomial ring and each face is identified with the product of the
corresponding variables.  A nonface is any subset of the vertices that does
not belong to the simplicial complex and each nonface is again identified with
a product of variables. The Stanley--Reisner ideal of a simplicial complex is
generated by the monomials corresponding to its nonfaces; see Definition~5.1.2
in \cite{BH}, Definition~1.6 in \cite{MS}, or Definition~II.1.1 in
\cite{Stanley}. Because computations in the associated polynomial ring are
typically prohibitive, this package is not intended for simplicial complexes
with a large number of vertices.

\subsection*{\scshape\mdseries Constructors}

\begin{figure}[t]
  \begin{subfigure}{0.3\textwidth}
    \begin{tikzpicture}[scale=0.7]
      \draw[fill=gray!50, ultra thick] (0,0) node[left]{\scriptsize $v$} --
      (0,1.50) node[left]{\scriptsize $w$} -- (0.75,0.75)
      node[below]{\scriptsize $x$} -- cycle;
      \draw[fill=gray!50, ultra thick] (1.50,0) node[right]{\scriptsize $y$} --
      (1.50,1.50) node[right]{\scriptsize $z$} -- (0.75,0.75) -- cycle;
      \filldraw (0,0) circle (2pt);
      \filldraw (0,1.50) circle (2pt);
      \filldraw (0.75,0.75) circle (2pt);
      \filldraw (1.50,0) circle (2pt);
      \filldraw (1.50,1.50) circle (2pt);
    \end{tikzpicture}
  \end{subfigure}
  \hspace{50pt}
  \begin{subfigure}{0.3\textwidth}
    \begin{tikzpicture}[scale=0.7]
      \draw[fill=gray!50, ultra thick] (0,0) -- (1.50,0) -- (1.50,1.50) --
      (0,1.50) -- cycle;
      \draw[ultra thick] (0,0) node[left]{\scriptsize $w$} -- (0,1.50)
      node[left]{\scriptsize $y$} -- (0.75,0.75)
      node[below]{\scriptsize $x$} -- cycle;
      \draw[fill=gray!50, ultra thick] (1.50,0) node[right]{\scriptsize $z$}
      -- (1.50,1.50) node[right]{\scriptsize $v$} -- (0.75,0.75) -- cycle;
      \filldraw (0,0) circle (2pt);
      \filldraw (0,1.50) circle (2pt);
      \filldraw (0.75,0.75) circle (2pt);
      \filldraw (1.50,0) circle (2pt);
      \filldraw (1.50,1.50) circle (2pt);
    \end{tikzpicture}
  \end{subfigure}
  \caption{On the left is the bowtie complex $⧓$ and on the right its Alexander dual
    $⧓^*$}
  \label{bowtie}
  \vspace{-1.0em}
\end{figure}
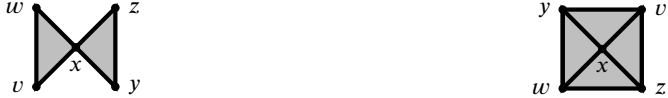

The basic constructor for a simplicial complex accepts two different kinds of
input.  Given a list of monomials, it returns the smallest simplicial complex
containing the corresponding faces.  Given a radical monomial ideal $I$, it
returns the simplicial complex whose Stanley-Reisner ideal is $I$.
We illustrate both methods using the `bowtie' complex in
Figure~\ref{bowtie}.

\begin{lstlisting}[xleftmargin=10pt, lineskip=-10pt, aboveskip=5.0pt, belowskip=1.5pt]
Macaulay2, version 1.20
with packages: ConwayPolynomials, Elimination, IntegralClosure, InverseSystems,
               Isomorphism, LLLBases, MinimalPrimes, OnlineLookup,
               PrimaryDecomposition, ReesAlgebra, Saturation, TangentCone
\end{lstlisting}
\begin{lstlisting}[xleftmargin=10pt, aboveskip=1.5pt, belowskip=1.5pt]
i1 : needsPackage "SimplicialComplexes"; S = QQ[v..z];
\end{lstlisting}
\begin{lstlisting}[xleftmargin=10pt, aboveskip=1.5pt, belowskip=1.5pt]
i3 : $⧓$ = simplicialComplex {v*w*x, x*y*z}
\end{lstlisting}
\begin{lstlisting}[xleftmargin=10pt, aboveskip=1.5pt, belowskip=1.5pt]
o3 = simplicialComplex | xyz vwx |
\end{lstlisting}
\begin{lstlisting}[xleftmargin=10pt, aboveskip=1.5pt, belowskip=1.5pt]
o3 : SimplicialComplex
\end{lstlisting}
\begin{lstlisting}[xleftmargin=10pt, aboveskip=1.5pt, belowskip=1.5pt]
i4 : I = monomialIdeal $⧓$
\end{lstlisting}
\begin{lstlisting}[xleftmargin=10pt, lineskip=-10pt, aboveskip=4pt, belowskip=1.5pt]
o4 = monomialIdeal (v*y, w*y, v*z, w*z)
\end{lstlisting}
\begin{lstlisting}[xleftmargin=10pt, aboveskip=1.5pt, belowskip=1.5pt]
o4 : MonomialIdeal of S
\end{lstlisting}
\begin{lstlisting}[xleftmargin=10pt, aboveskip=1.5pt, belowskip=1.5pt]
i5 : $⧑$ = simplicialComplex I
\end{lstlisting}
\begin{lstlisting}[xleftmargin=10pt, aboveskip=1.5pt, belowskip=1.5pt]
o5 = simplicialComplex | xyz vwx |
\end{lstlisting}
\begin{lstlisting}[xleftmargin=10pt, aboveskip=1.5pt, belowskip=1.5pt]
o5 : SimplicialComplex
\end{lstlisting}
\begin{lstlisting}[xleftmargin=10pt, aboveskip=1.5pt, belowskip=3.0pt]
i6 : assert($⧓$ === $⧑$)
\end{lstlisting}

The package also has convenient constructors for some archetypal simplicial
complexes. For example, we recognize the real projective plane
and the Klein bottle from the reduced homology groups of some explicit
triangulations; see Theorems~6.3--6.4 in \cite{Munkres}.
\begin{lstlisting}[xleftmargin=10pt, aboveskip=3.0pt, belowskip=1.5pt]
i7 : $ℙ$ = realProjectiveSpaceComplex(2, R = ZZ[a..h])
\end{lstlisting}
\begin{lstlisting}[xleftmargin=10pt, aboveskip=1.5pt, belowskip=1.5pt]
o7 = simplicialComplex | bef aef cdf adf bcf cde bde ace abd abc |
\end{lstlisting}
\begin{lstlisting}[xleftmargin=10pt, aboveskip=1.5pt, belowskip=1.5pt]
o7 : SimplicialComplex
\end{lstlisting}
\begin{lstlisting}[xleftmargin=10pt, aboveskip=1.5pt, belowskip=1.5pt]
i8 : for j from 0 to 2 list prune HH_j $ℙ$
\end{lstlisting}
\begin{lstlisting}[xleftmargin=10pt, aboveskip=1.5pt, belowskip=1.5pt]
o8 = {0, cokernel | 2 |, 0}
\end{lstlisting}
\begin{lstlisting}[xleftmargin=10pt, aboveskip=1.5pt, belowskip=1.5pt]
o8 : List
\end{lstlisting}
\begin{lstlisting}[xleftmargin=10pt, aboveskip=1.5pt, belowskip=1.5pt]
i9 : for j from 0 to 2 list prune HH_j kleinBottleComplex R
\end{lstlisting}
\begin{lstlisting}[xleftmargin=10pt, lineskip=-10pt, aboveskip=4pt, belowskip=1.5pt]
o9 = {0, cokernel | 2 |, 0}
                  | 0 |
\end{lstlisting}
\begin{lstlisting}[xleftmargin=10pt, aboveskip=-5pt, belowskip=3.0pt]
o9 : List
\end{lstlisting}
More comprehensively, Frank H.~Lutz enumerates simplicial complexes having a
small number of vertices; see \cite{LutzM}.  Using this list, the package
creates a database of $43\,138$ simplicial $2$\nobreakdash-manifolds having at
most $10$ vertices and $1\,343$ simplicial $3$-manifolds having at most $9$
vertices.  We demonstrate this feature by exhibiting the distribution of
$f\!$-vectors among the $3$-manifolds having $9$ vertices.  For all nonnegative
integers $j$, the $j$-th entry in the $f\!$-vector is the number of faces having
$j$ vertices.
\begin{lstlisting}[xleftmargin=10pt, aboveskip=5.0pt, belowskip=1.5pt]
i10 : tally for j from 0 to 1296 list fVector smallManifold(3, 9, j, ZZ[vars(1..9)])
\end{lstlisting}
\begin{lstlisting}[xleftmargin=10pt, lineskip=-10pt, aboveskip=4pt, belowskip=1.5pt]
o10 = Tally{{1, 9, 26, 34, 17} => 7  }
            {1, 9, 27, 36, 18} => 23
            {1, 9, 28, 38, 19} => 45
            {1, 9, 29, 40, 20} => 84
            {1, 9, 30, 42, 21} => 128
            {1, 9, 31, 44, 22} => 175
            {1, 9, 32, 46, 23} => 223
            {1, 9, 33, 48, 24} => 231
            {1, 9, 34, 50, 25} => 209
            {1, 9, 35, 52, 26} => 121
            {1, 9, 36, 54, 27} => 51
\end{lstlisting}
\begin{lstlisting}[xleftmargin=10pt, aboveskip=-5pt, belowskip=3.0pt]
o10 : Tally
\end{lstlisting}
Exploiting the same loop, we construct the simplicial maps from a minimal
triangulation of a torus to the induced subcomplex on the first $7$ vertices
for each of these $3$-manifolds.
\begin{lstlisting}[xleftmargin=10pt, aboveskip=3.0pt, belowskip=1.5pt]
i11 : T = smallManifold(2, 7, 6, R = ZZ[a..i])
\end{lstlisting}
\begin{lstlisting}[xleftmargin=10pt, aboveskip=1.5pt, belowskip=1.5pt]
o11 = simplicialComplex | cfg afg beg aeg cdg bdg def bef adf bcf cde ace abd abc |
\end{lstlisting}
\begin{lstlisting}[xleftmargin=10pt, aboveskip=1.5pt, belowskip=1.5pt]
o11 : SimplicialComplex
\end{lstlisting}
\begin{lstlisting}[xleftmargin=10pt, aboveskip=1.5pt, belowskip=1.5pt]
i12 : for j from 0 to 2 list prune HH_j T
\end{lstlisting}
\begin{lstlisting}[xleftmargin=10pt, lineskip=-10pt, aboveskip=4pt, belowskip=1.5pt]
            2   1
o12 = {0, ZZ, ZZ }
\end{lstlisting}
\begin{lstlisting}[xleftmargin=10pt, aboveskip=1.5pt, belowskip=1.5pt]
o12 : List
\end{lstlisting}
\begin{lstlisting}[xleftmargin=10pt, lineskip=-10pt, aboveskip=4pt, belowskip=1.5pt]
i13 : for j from 0 to 1296 list (
          phi := map(smallManifold(3, 9, j, R), T, gens R);
          if not isWellDefined phi then continue else phi);
\end{lstlisting}
\begin{lstlisting}[xleftmargin=10pt, aboveskip=1.5pt, belowskip=1.5pt]
o13 : {}
\end{lstlisting}
\begin{lstlisting}[xleftmargin=10pt, aboveskip=1.5pt, belowskip=-1.0pt]
o13 : List
\end{lstlisting}

\subsection*{\scshape\mdseries Combinatorial Topology}

We use the bowtie complex to showcase some of the key operations on simplicial
complexes. Viewing a simplicial complex as a subcomplex of a simplex yields a
duality theory. For any simplicial complex $\Delta$ whose vertices belong to a
set $V$, the Alexander dual is the simplicial complex
$\Delta^* \colequal \{ F \subseteq V \mathrel{|} V \setminus F \not\in \Delta \}$.
Since each simplicial complex in this package has an underlying polynomial
ring, the variables in this ring form a canonical superset of the vertices.
\begin{lstlisting}[xleftmargin=10pt, aboveskip=3.0pt, belowskip=1.5pt]
i14 : dual $⧓$
\end{lstlisting}
\begin{lstlisting}[xleftmargin=10pt, aboveskip=1.5pt, belowskip=1.5pt]
o14 = simplicialComplex | wxz vxz wxy vxy |
\end{lstlisting}
\begin{lstlisting}[xleftmargin=10pt, aboveskip=1.5pt, belowskip=1.5pt]
o14 : SimplicialComplex
\end{lstlisting}
\begin{lstlisting}[xleftmargin=10pt, aboveskip=1.5pt, belowskip=1.5pt]
i15 : assert(dual dual $⧓$ === $⧓$ and dual monomialIdeal $⧓$ === monomialIdeal dual $⧓$)
\end{lstlisting}
Algebraically, Alexander duality switches the roles of the minimal generators
and the irreducible components in the Stanley--Reisner ideal.
\begin{lstlisting}[xleftmargin=10pt, aboveskip=3pt, belowskip=1.5pt]
i16 : monomialIdeal dual $⧓$
\end{lstlisting}
\begin{lstlisting}[xleftmargin=10pt, aboveskip=4pt, belowskip=1pt]
o16 = monomialIdeal (v*w, y*z)
\end{lstlisting}
\begin{lstlisting}[xleftmargin=10pt, aboveskip=1.5pt, belowskip=1.5pt]
o16 : MonomialIdeal of S
\end{lstlisting}
\begin{lstlisting}[xleftmargin=10pt, aboveskip=1.5pt, belowskip=1.5pt]
i17 : irreducibleDecomposition monomialIdeal $⧓$
\end{lstlisting}
\begin{lstlisting}[xleftmargin=10pt, aboveskip=4pt, belowskip=1pt]
o17 = {monomialIdeal (v, w), monomialIdeal (y, z)}
\end{lstlisting}
\begin{lstlisting}[xleftmargin=10pt, aboveskip=1pt, belowskip=3pt]
o17 : List
\end{lstlisting}
The topological form of Alexander duality gives an isomorphism
between the reduced homology of a simplicial complex and reduced cohomology of
its dual; see Theorem~5.6 in \cite{MS}.
\begin{lstlisting}[xleftmargin=10pt, aboveskip=3pt, belowskip=1.5pt]
i18 : n = numgens ring $⧓$
\end{lstlisting}
\begin{lstlisting}[xleftmargin=10pt, aboveskip=4pt, belowskip=1pt]
o18 = 5
\end{lstlisting}
\begin{lstlisting}[xleftmargin=10pt, aboveskip=1pt, belowskip=3pt]
i19 : assert all(-1..n-1, j -> prune HH^(n-j-3) dual $⧓$ == prune HH_j $⧓$)
\end{lstlisting}

A simplicial complex $\Delta$ is Cohen--Macaulay if the associated quotient
ring $S/I$, where $I$ is the Stanley--Reisner ideal of $\Delta$ in the
polynomial ring $S$, is Cohen--Macaulay.  To characterize this attribute
topologically, we introduce a family of subcomplexes.  For any face $F$ in
$\Delta$, the link is the subcomplex
$\operatorname{link}_\Delta(F) \colequal \{ G \in \Delta \mathrel{|}
\text{$F \cup G \in \Delta$ and $F \cap G = \varnothing$} \}$.  The link of
the vertex $x$ in $⧓$ has two disjoint facets.
\begin{lstlisting}[xleftmargin=10pt, aboveskip=3.0pt, belowskip=1.5pt]
i20 : L = link($⧓$, x)
\end{lstlisting}
\begin{lstlisting}[xleftmargin=10pt, aboveskip=1.5pt, belowskip=1.5pt]
o20 = simplicialComplex | yz vw |
\end{lstlisting}
\begin{lstlisting}[xleftmargin=10pt, aboveskip=1.5pt, belowskip=1.5pt]
o20 : SimplicialComplex
\end{lstlisting}
\begin{lstlisting}[xleftmargin=10pt, aboveskip=1.5pt, belowskip=1.5pt]
i21 : prune HH_0 L
\end{lstlisting}
\begin{lstlisting}[xleftmargin=10pt, lineskip=-10pt, aboveskip=4pt, belowskip=1pt]
        1
o21 = QQ
\end{lstlisting}
\begin{lstlisting}[xleftmargin=10pt, aboveskip=1.5pt, belowskip=3.0pt]
o21 : QQ-module, free
\end{lstlisting}
As discovered by Gerald Reisner, the simplicial complex $\Delta$ is
Cohen--Macaulay if and only if, for all faces $F$ in $\Delta$ and all integers
$j$ less than the dimension of $\operatorname{link}_\Delta(F)$, the $j$-th
reduced homology group of $\operatorname{link}_\Delta(F)$ vanishes; see
Corollary~5.3.9 in \cite{BH}, Theorem~5.53 in \cite{MS}, or Corollary~II.4.2
in \cite{Stanley}.  Using this criterion, the $0$-th reduced homology
certifies that $⧓$ is not Cohen--Macaulay.
\begin{lstlisting}[xleftmargin=10pt, aboveskip=3.0pt, belowskip=1.5pt]
i22 : assert(HH_0 L != 0)
\end{lstlisting}
\begin{lstlisting}[xleftmargin=10pt, aboveskip=1.5pt, belowskip=3.0pt]
i23 : assert(dim(S^1/monomialIdeal $⧓$) =!= n - pdim(S^1/monomialIdeal $⧓$))
\end{lstlisting}
However, the $1$-skeleton of $⧓$ is Cohen--Macaulay.
\begin{lstlisting}[xleftmargin=10pt, aboveskip=3.0pt, belowskip=1.5pt]
i24 : $⋈$ = skeleton(1, $⧓$)
\end{lstlisting}
\begin{lstlisting}[xleftmargin=10pt, aboveskip=1.5pt, belowskip=1.5pt]
o24 = simplicialComplex | yz xz xy wx vx vw |
\end{lstlisting}
\begin{lstlisting}[xleftmargin=10pt, aboveskip=1.5pt, belowskip=1.5pt]
o24 : SimplicialComplex
\end{lstlisting}
\begin{lstlisting}[xleftmargin=10pt, aboveskip=1.5pt, belowskip=1.5pt]
i25 : faceList = rsort flatten values faces $⋈$
\end{lstlisting}
\begin{lstlisting}[xleftmargin=10pt, aboveskip=4pt, belowskip=1pt]
o25 = {v*w, v*x, w*x, x*y, x*z, y*z, v, w, x, y, z, 1}
\end{lstlisting}
\begin{lstlisting}[xleftmargin=10pt, aboveskip=1.5pt, belowskip=1.5pt]
o25 : List
\end{lstlisting}
\begin{lstlisting}[xleftmargin=10pt, aboveskip=1.5pt, belowskip=1.5pt]
i26 : assert all(faceList, F -> (L := link($⋈$, F); all(dim L, j -> HH_j L == 0)))
\end{lstlisting}
\begin{lstlisting}[xleftmargin=10pt, aboveskip=1.5pt, belowskip=3.0pt]
i27 : assert(dim(S^1/monomialIdeal $⋈$) === n - pdim(S^1/monomialIdeal $⋈$))
\end{lstlisting}

Alternatively, we verify that $⧓$ is not Cohen--Macaulay by showing that its
$h$-vector has a negative entry; see Theorem~5.1.10 in \cite{BH} or
Corollary~II.2.5 in \cite{Stanley}.  By definition, the $h$-vector of a
simplicial complex $\Delta$ is a binomial transform of its $f\!$-vector: for
all $0 \leqslant j \leqslant d \colequal \dim \Delta$, we have
$h_j = \sum_{k=0}^{j} (-1)^{j-1}
\smash{\tbinom{d+1-k}{\raisebox{3pt}{$\scriptstyle j-k$}}} f_{k-1}$.  The
$h$-vector encodes the numerator of the Hilbert series for $S/I$.
\begin{lstlisting}[xleftmargin=10pt, aboveskip=5.0pt, belowskip=1.5pt]
i28 : d = dim $⧓$
\end{lstlisting}
\begin{lstlisting}[xleftmargin=10pt, aboveskip=1.5pt, belowskip=1.5pt]
o28 = 2
\end{lstlisting}
\begin{lstlisting}[xleftmargin=10pt, aboveskip=1.5pt, belowskip=1.5pt]
i29 : faces $⧓$
\end{lstlisting}
\begin{lstlisting}[xleftmargin=10pt, lineskip=-10pt, aboveskip=4pt, belowskip=1pt]
o29 =  HashTable{-1 => {1}                          }
                0 => {v, w, x, y, z}
                1 => {v*w, v*x, w*x, x*y, x*z, y*z}
                2 => {v*w*x, x*y*z}
\end{lstlisting}
\begin{lstlisting}[xleftmargin=10pt, aboveskip=1.5pt, belowskip=1.5pt]
o29 : HashTable
\end{lstlisting}
\begin{lstlisting}[xleftmargin=10pt, aboveskip=1.5pt, belowskip=1.5pt]
i30 : fVec = fVector $⧓$
\end{lstlisting}
\begin{lstlisting}[xleftmargin=10pt, aboveskip=1.5pt, belowskip=1.5pt]
o30 = {1, 5, 6, 2}
\end{lstlisting}
\begin{lstlisting}[xleftmargin=10pt, aboveskip=1.5pt, belowskip=1.5pt]
o30 : List
\end{lstlisting}
\begin{lstlisting}[xleftmargin=10pt, lineskip=-10pt, aboveskip=4pt, belowskip=1pt]
i31 : hVec = for j from 0 to d list 
          sum(j+1, k -> (-1)^(j-k) * binomial(d+1-k, j-k) * fVec#k)
\end{lstlisting}
\vfill
\begin{lstlisting}[xleftmargin=10pt, aboveskip=1.5pt, belowskip=1.5pt]
o31 = {1, 2, -1}
\end{lstlisting}
\begin{lstlisting}[xleftmargin=10pt, aboveskip=1.5pt, belowskip=1.5pt]
o31 : List
\end{lstlisting}
\begin{lstlisting}[xleftmargin=10pt, aboveskip=1.5pt, belowskip=1.5pt]
i32 : hilbertSeries(S^1/monomialIdeal $⧓$, Reduce => true)
\end{lstlisting}
\begin{lstlisting}[xleftmargin=10pt, lineskip=-10pt, aboveskip=4pt, belowskip=1pt]
                2
      1 + 2T - T
o32 = -----------
               3
        (1 - T)
\end{lstlisting}
\begin{lstlisting}[xleftmargin=10pt, aboveskip=1.5pt, belowskip=-5.0pt]
o32 : Expression of class Divide
\end{lstlisting}

\subsection*{\scshape\mdseries Resolutions of Monomial Ideals}

As David Bayer, Irena Peeva, and Bernd Sturmfels~\cite{BPS} reveal, minimal
free resolutions of monomial ideals are frequently encoded by a simplicial
complex.  Consider a monomial ideal $J$ in the polynomial ring
$R \colequal \mathbb{Q}[y_1, y_2, \dotsc, y_m]$.  Assume that $R$ is equipped
with the $\mathbb{N}^{m}$-grading given by $\deg(y_i) = \textbf{e}_i$, for each
$1 \leqslant i \leqslant m$, where
$\textbf{e}_1, \textbf{e}_2, \dotsc, \textbf{e}_m$ is the standard basis.
Let $\Delta$ be a simplicial complex whose vertices are labelled by the generators
of $J$.  We label each face $F$ of $\Delta$ by the least common multiple
$y^{\,\textbf{a}_{\!F}} \in R$ of its vertices; the empty face is labelled by
the monomial $1 = y^{\,\textbf{a}_{\!\varnothing}}$.  The chain complex
$C(\Delta)$ supported on the labelled simplicial complex $\Delta$ is the chain
complex of free $\mathbb{N}^m$-graded $R$-modules with basis corresponding to
the faces of $\Delta$. More precisely, the chain complex $C(\Delta)$ is determined
by the data
\begin{align*}
  C_i(\Delta)
  &\colequal \bigoplus_{\dim(F) \,=\, i-1} R(-\textbf{a}_F) 
  & & \text{and}
  &\partial(F)
  &= \sum_{\dim(G) \,=\, \dim(F)-1} \operatorname{sign}(G,F) \,
    y^{\,\textbf{a}_{\!F} - \textbf{a}_{\!G}} \, G \, .
\end{align*}
The symbols $F$ and $G$ represent both faces in $\Delta$ and basis vectors in
the underlying free module of $C(\Delta)$.  The sign of the pair $(G,F)$
belongs to $\{-1,0,1\}$ and is part of the data in the boundary map of the
chain complex of $\Delta$.  For more information, see Subsection~4.1 in
\cite{MS} or Chapter~55 in \cite{Peeva}.

\begin{figure}[!ht]
  \centering
  \begin{subfigure}{0.3\textwidth}
    \begin{tikzpicture}
      \filldraw[gray!50] (0,0) -- (0,-1) -- (0.866025,-0.5) -- cycle;
      \draw[ultra thick] (0,0) node[left]{\scriptsize $x_0$} -- (0,-1)
      node[left]{\scriptsize $x_1$} -- (0.866025,-0.5)
      node[below]{\scriptsize $\,\,\,x_2$} -- cycle;
      \draw[ultra thick] (0.866025,-0.5) -- (1.866025,-0.5)
      node[right]{\scriptsize $x_3$};
      \filldraw (0,0) circle (2pt);
      \filldraw (0,-1) circle (2pt);
      \filldraw (0.866025,-0.5) circle (2pt);
      \filldraw (1.866025,-0.5) circle (2pt);
    \end{tikzpicture}
  \end{subfigure}
  \hspace{50pt}
  \begin{subfigure}{0.3\textwidth}
    \begin{tikzpicture}
      \filldraw[gray!50] (0,0) -- (0,-1) -- (0.866025,-0.5) -- cycle;
      \draw[ultra thick] (0,0) node[left]{\scriptsize $y_0 y_1$} -- (0,-1)
      node[left]{\scriptsize $y_0 y_2$} -- (0.866025,-0.5)
      node[below]{\scriptsize $\,\,\,y_0 y_3$} -- cycle;
      \draw[ultra thick] (0.866025,-0.5) -- (1.866025,-0.5)
      node[right]{\scriptsize $y_1 y_2 y_3$};
      \filldraw (0,0) circle (2pt);
      \filldraw (0,-1) circle (2pt);
      \filldraw (0.866025,-0.5) circle (2pt);
      \filldraw (1.866025,-0.5) circle (2pt);
    \end{tikzpicture}
  \end{subfigure}
  \caption{On the left is $\Gamma$ and on the right is the
    labelling of its vertices}
  \label{Gamma}
\end{figure}
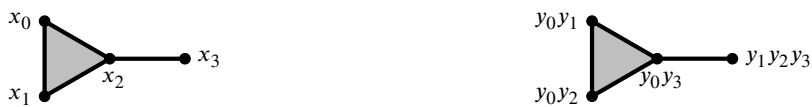

We illustrate this construction with an explicit example.  Consider the
simplicial complex $\Gamma$ in Figure~\ref{Gamma} and the monomial ideal
$J = (y_0 y_1, y_0 y_2, y_0 y_3, y_1 y_2 y_3)$ in
$R = \mathbb{Q}[y_0,y_1, y_2, y_3]$.
Label the vertices of $\Gamma$ by the generators of $J$:
$x_0 \mapsto y_0 y_1$, $x_1 \mapsto y_0 y_2$, $x_2 \mapsto y_0 y_3$, and
$x_3 \mapsto y_1 y_2 y_3$.  
\begin{lstlisting}[xleftmargin=10pt, aboveskip=3.0pt, belowskip=1.5pt]
i33 : x = getSymbol "x"; S = ZZ[x_0..x_3];
\end{lstlisting}
\begin{lstlisting}[xleftmargin=10pt, aboveskip=1.5pt, belowskip=1.5pt]
i35 : $Δ$ = simplicialComplex{x_0*x_1*x_2, x_2*x_3}
\end{lstlisting}
\begin{lstlisting}[xleftmargin=10pt, aboveskip=1.5pt, belowskip=1.5pt]
o35 = simplicialComplex | x_2x_3 x_0x_1x_2 |
\end{lstlisting}
\begin{lstlisting}[xleftmargin=10pt, aboveskip=1.5pt, belowskip=1.5pt]
o35 : SimplicialComplex
\end{lstlisting}
\begin{lstlisting}[xleftmargin=10pt, aboveskip=1.5pt, belowskip=1.5pt]
i36 : chainComplex $Δ$
\end{lstlisting}
\begin{lstlisting}[xleftmargin=10pt, lineskip=-10pt, aboveskip=1.5pt, belowskip=1.5pt]
        1       4       4       1
o36 = ZZ  <-- ZZ  <-- ZZ  <-- ZZ
\end{lstlisting}
\begin{lstlisting}[xleftmargin=10pt, aboveskip=1.5pt, belowskip=1.5pt]
      -1      0       1       2
\end{lstlisting}
\begin{lstlisting}[xleftmargin=10pt, aboveskip=1.5pt, belowskip=1.5pt]
o36 : ChainComplex
\end{lstlisting}
\begin{lstlisting}[xleftmargin=10pt, aboveskip=1.5pt, belowskip=1.5pt]
i37 : y = getSymbol "y"; R = QQ[y_0..y_3, DegreeRank => 4];
\end{lstlisting}
\begin{lstlisting}[xleftmargin=10pt, aboveskip=1.5pt, belowskip=1.5pt]
i39 : J = ideal(y_0*y_1, y_0*y_2, y_0*y_3, y_1*y_2*y_3)
\end{lstlisting}
\begin{lstlisting}[xleftmargin=10pt, lineskip=-10pt, aboveskip=1.5pt, belowskip=1.5pt]
o39 = ideal (y y , y y , y y , y y y )
              0 1   0 2   0 3   1 2 3
\end{lstlisting}
\begin{lstlisting}[xleftmargin=10pt, aboveskip=1.5pt, belowskip=1.5pt]
o39 : Ideal of R
\end{lstlisting}
\begin{lstlisting}[xleftmargin=10pt, aboveskip=1.5pt, belowskip=1.5pt]
i40 : C = chainComplex($Δ$, Labels => J_*)
\end{lstlisting}
\begin{lstlisting}[xleftmargin=10pt, lineskip=-10pt, aboveskip=1.5pt, belowskip=1.5pt]
       1      4      4      1
o40 = R  <-- R  <-- R  <-- R
\end{lstlisting}
\begin{lstlisting}[xleftmargin=10pt, aboveskip=1.5pt, belowskip=1.5pt]                           
      0      1      2      3
\end{lstlisting}
\begin{lstlisting}[xleftmargin=10pt, aboveskip=1.5pt, belowskip=1.5pt]
o40 : ChainComplex
\end{lstlisting}
\begin{lstlisting}[xleftmargin=10pt, aboveskip=1.5pt, belowskip=1.5pt]
i41 : C.dd
\end{lstlisting}
\begin{lstlisting}[xleftmargin=10pt, lineskip=-10pt, aboveskip=1.5pt, belowskip=3.0pt]
           1                                          4
o41 = 0 : R  <-------------------------------------- R  : 1
               | y_0y_1 y_0y_2 y_0y_3 y_1y_2y_3 |
\end{lstlisting}
\begin{lstlisting}[xleftmargin=10pt, lineskip=-10pt, aboveskip=1.5pt, belowskip=3.0pt]
           4                                               4
      1 : R  <------------------------------------------- R  : 2
                {1, 1, 0, 0} | -y_2 -y_3 0    0       |
                {1, 0, 1, 0} | y_1  0    -y_3 0       |
                {1, 0, 0, 1} | 0    y_1  y_2  -y_1y_2 |
                {0, 1, 1, 1} | 0    0    0    y_0     |
\end{lstlisting}
\begin{lstlisting}[xleftmargin=10pt, lineskip=-10pt, aboveskip=1.5pt, belowskip=1.5pt]
           4                             1
      2 : R  <------------------------- R  : 3
                {1, 1, 1, 0} | y_3  |
                {1, 1, 0, 1} | -y_2 |
                {1, 0, 1, 1} | y_1  |
                {1, 1, 1, 1} | 0    |
\end{lstlisting}
\begin{lstlisting}[xleftmargin=10pt, aboveskip=1.5pt, belowskip=1.5pt]
o41 : ChainComplexMap
\end{lstlisting}
\begin{lstlisting}[xleftmargin=10pt, aboveskip=1.5pt, belowskip=3.0pt]
i42 : assert(res(R^1/J) == C)
\end{lstlisting}
The chain complex $C(\Delta)$ 
depends on the labelling and is not always a resolution.
\begin{lstlisting}[xleftmargin=10pt, aboveskip=3.0pt, belowskip=1.5pt]
i43 : C' = chainComplex($Δ$, Labels => reverse J_*)
\end{lstlisting}
\begin{lstlisting}[xleftmargin=10pt, lineskip=-10pt, aboveskip=1.5pt, belowskip=1.5pt]
       1      4      4      1
o43 = R  <-- R  <-- R  <-- R
\end{lstlisting}
\begin{lstlisting}[xleftmargin=10pt, aboveskip=1.5pt, belowskip=1.5pt]
      0      1      2      3
\end{lstlisting}
\begin{lstlisting}[xleftmargin=10pt, aboveskip=1.5pt, belowskip=1.5pt] 
o43 : ChainComplex
\end{lstlisting}
\begin{lstlisting}[xleftmargin=10pt, aboveskip=1.5pt, belowskip=1.5pt] 
i44 : prune homology C'
\end{lstlisting}
\begin{lstlisting}[xleftmargin=10pt, aboveskip=1.5pt, belowskip=1.5pt]
o44 = 0 : cokernel | y_0y_3 y_0y_2 y_0y_1 y_1y_2y_3 |
\end{lstlisting}
\begin{lstlisting}[xleftmargin=10pt, aboveskip=1.5pt, belowskip=1.5pt]
      1 : cokernel {1, 1, 0, 1} | y_2 |              
\end{lstlisting}
\begin{lstlisting}[xleftmargin=10pt, aboveskip=1.5pt, belowskip=1.5pt]
      2 : 0                                          
\end{lstlisting}
\begin{lstlisting}[xleftmargin=10pt, aboveskip=1.5pt, belowskip=1.5pt]
      3 : 0   
\end{lstlisting}
\begin{lstlisting}[xleftmargin=10pt, aboveskip=1.5pt, belowskip=3.0pt]
o44 : GradedModule
\end{lstlisting}

Given a monomial ideal $J$, there are several algorithms that return a
labelled simplicial complex $\Delta$ such that chain complex $C(\Delta)$ is a
free resolution of $R/J$.  We exhibit a few.
\begin{lstlisting}[xleftmargin=10pt, aboveskip=3.0pt, belowskip=1.5pt]
i45 : J' = monomialIdeal(y_1*y_3, y_2^2, y_0*y_2, y_1^2, y_0^2);
\end{lstlisting}
\begin{lstlisting}[xleftmargin=10pt, aboveskip=1.5pt, belowskip=1.5pt]
o45 : MonomialIdeal of R
\end{lstlisting}
\begin{lstlisting}[xleftmargin=10pt, aboveskip=1.5pt, belowskip=1.5pt]
i46 : T = taylorResolution J'
\end{lstlisting}
\vfill
\begin{lstlisting}[xleftmargin=10pt, lineskip=-10pt, aboveskip=1.5pt, belowskip=1.5pt]
       1      5      10      10      5      1
o46 = R  <-- R  <-- R   <-- R   <-- R  <-- R
\end{lstlisting}
\begin{lstlisting}[xleftmargin=10pt, aboveskip=1.5pt, belowskip=1.5pt] 
      0      1      2       3       4      5
\end{lstlisting}
\begin{lstlisting}[xleftmargin=10pt, aboveskip=1.5pt, belowskip=1.5pt]
o46 : ChainComplex
\end{lstlisting}
\begin{lstlisting}[xleftmargin=10pt, aboveskip=1.5pt, belowskip=1.5pt]
i47 : gensJ' = first entries mingens J'
\end{lstlisting}
\begin{lstlisting}[xleftmargin=10pt, lineskip=-10pt, aboveskip=1.5pt, belowskip=1.5pt]
              2         2   2
o47 = {y y , y , y y , y , y }
        1 3   2   0 2   1   0
\end{lstlisting}
\begin{lstlisting}[xleftmargin=10pt, aboveskip=1.5pt, belowskip=1.5pt]
o47 : List
\end{lstlisting}
\begin{lstlisting}[xleftmargin=10pt, aboveskip=1.5pt, belowskip=1.5pt]
i48 : S = ZZ[x_0..x_4];
\end{lstlisting}
\begin{lstlisting}[xleftmargin=10pt, aboveskip=1.5pt, belowskip=1.5pt]
i49 : assert(T == chainComplex(simplexComplex(4, S), Labels => gensJ'))
\end{lstlisting}
\begin{lstlisting}[xleftmargin=10pt, aboveskip=1.5pt, belowskip=1.5pt]
i50 : assert(lyubeznikSimplicialComplex(J', S) === simplexComplex(4, S))
\end{lstlisting}
\begin{lstlisting}[xleftmargin=10pt, aboveskip=1.5pt, belowskip=1.5pt]
i51 : $Γ$ = buchbergerSimplicialComplex(J',S)
\end{lstlisting}
\begin{lstlisting}[xleftmargin=10pt, aboveskip=1.5pt, belowskip=1.5pt]
o52 = simplicialComplex | x_0x_2x_3x_4 x_0x_1x_2x_3 |
\end{lstlisting}
\begin{lstlisting}[xleftmargin=10pt, aboveskip=1.5pt, belowskip=1.5pt]
o52 : SimplicialComplex
\end{lstlisting}
\begin{lstlisting}[xleftmargin=10pt, aboveskip=1.5pt, belowskip=1.5pt]
i53 : B = buchbergerResolution J'
\end{lstlisting}
\begin{lstlisting}[xleftmargin=10pt, lineskip=-10pt, aboveskip=1.5pt, belowskip=1.5pt]
       1      5      9      7      2
o53 = R  <-- R  <-- R  <-- R  <-- R
\end{lstlisting}
\begin{lstlisting}[xleftmargin=10pt, aboveskip=1.5pt, belowskip=1.5pt]
      0      1      2      3      4
\end{lstlisting}
\begin{lstlisting}[xleftmargin=10pt, aboveskip=1.5pt, belowskip=1.5pt]
o53 : ChainComplex
\end{lstlisting}
\begin{lstlisting}[xleftmargin=10pt, aboveskip=1.5pt, belowskip=1.5pt]
i54 : assert all(3, i -> HH_(i+1) B == 0) 
\end{lstlisting}
\begin{lstlisting}[xleftmargin=10pt, aboveskip=1.5pt, belowskip=1.5pt]
i55 : assert(betti B == betti res J')
\end{lstlisting}
\begin{lstlisting}[xleftmargin=10pt, aboveskip=1.5pt, belowskip=1.5pt]
i56 : assert(B == chainComplex($Γ$, Labels => first entries mingens J'))
\end{lstlisting}
\begin{lstlisting}[xleftmargin=10pt, aboveskip=1.5pt, belowskip=1.5pt]
i57 : assert($Γ$ === lyubeznikSimplicialComplex(J', S, MonomialOrder => {2,1,0,3,4}))
\end{lstlisting}
\begin{lstlisting}[xleftmargin=10pt, aboveskip=1.5pt, belowskip=3.0pt]
i59 : assert($Γ$ === scarfSimplicialComplex(J', S))
\end{lstlisting}
For more information about the Taylor resolution, the Lyubeznik resolution,
and the Scarf complex, see \cite{Mermin}.  The Buchberger resolution is
described in \cite{OW}.  

\subsection*{\scshape\mdseries Acknowledgements}
All three authors were partially supported by the Natural Sciences and
Engineering Research Council of Canada (NSERC).


\vspace*{-1em}
\raggedright

\end{document}